# An elementary algorithm for classifying digital objects based on the variational principle.

M.A. Antonets

## 1. Introduction.

In [1-2], it was proposed to use the maximin principle to solve the problem of classifying various digital objects by histograms generated by these objects. Examples of such objects are texts, raster images, cardiograms, and other numerical dependencies studied in computer science, medicine, and engineering. As in game theory, the solution to the optimal course of action problem is given by a weight function on a finite set [3]. That weight functions are used to define weight of any histogram is considered as a quantitative measure of the correspondence of an object to a given class, and the maximin value is the threshold for the acceptable weight of the histogram. This weight is optimal because it leads to the most stringent rule for selecting objects using weighting. However, one should take into account the possibility of having several vertices of a polyhedron generated by inequalities of constraints on which the values of the objective function are equal to or close to the value of maximin. This may serve as a basis for additional analysis of the expertise results (see [2]).

The paper considers a linear programming problem in the following form: for given non-negative numbers $b_1, \dots, b_m$ real functions $c, a_1, \dots, a_m$ on the set $\Omega_1 \stackrel{\text{def}}{=} \{1, \dots, n+1\}$ it is necessary to find all non-negative functions $x$ on the set $\Omega_1$ for which the sum
$c(1)x(1) + \cdots + c(n)x(n)$ reaches a maximum subject to the inequalities-constraints
$a_l(1)x(1) + \cdots + a_l(n)x(n) \leq b_k$, $l = 1, \dots, m$.
Algorithms for finding all solutions to this problem were proposed in the works [4]-[8].
H. Uzawa obtained in [4] explicit expressions for all solutions of the problem under consideration, containing products of matrices whose sizes increase with increasing dimension of the constraint vectors $n$ and their number $m$.

## 2. Formulas for solutions of the linear programming problems.

Following work [4], we transform the original data by introducing the functions $f^l, l = 1, \dots, m-1$ on the set $\Omega_1$ defined by relation

$$f^l \stackrel{\text{def}}{=} \left( a_{l+1}(1) - \frac{b(l)}{b(1)} a_1(1), \dots, a_{l+1}(n) - \frac{b(l)}{b(1)} a_1(n), -\frac{b(l)}{b(1)} \right) \quad (1)$$

Where $b_1 \neq 0$.

Having the goal of constructing sequential sets $\Omega_l$ for $l \in [2, \dots, m-1]$ and functions $g^l$ on each of these sets we put $g^1 \stackrel{\text{def}}{=} f^1$ and define sets $\Omega_1^0, \Omega_1^-, \Omega_1^+$ by relations

$$\Omega_1^0 \stackrel{\text{def}}{=} \{k \in \Omega_1, g^1(k) = 0\}, \quad \Omega_1^- \stackrel{\text{def}}{=} \{k \in \Omega_1, g^1(k) < 0\}, \quad \Omega_1^+ \stackrel{\text{def}}{=} \{k \in \Omega_1, g^1(k) > 0\}$$

and set $\Omega_2 \stackrel{\text{def}}{=} \widehat{\Omega_1^0 \cup \Omega_1^-} \cup \Omega_1^{\mp}$, where $\widehat{\Omega_1^0 \cup \Omega_1^-} \stackrel{\text{def}}{=} \{1,2,\dots,|\Omega_1^0 \cup \Omega_1^-|\}$, $\Omega_1^{\mp} \stackrel{\text{def}}{=} \Omega_1^- \times \Omega_1^+$.

Finally, let's denote $j_1(\omega_1) \stackrel{\text{def}}{=} \omega_1$, $\omega_1 \in \Omega_1$ and define the numbering $j_2$ on the set $\Omega_2$ as follows: set $j_2(\omega_2) \stackrel{\text{def}}{=} \omega_2$ for $\omega_2 \in \widehat{\Omega_1^0 \cup \Omega_1^-}$. Then, considering all elements from $\Omega_1^{\mp}$ to follow elements from the set $\Omega_1^0 \cup \Omega_1^-$, we will assume that of two pairs $\{\omega_-', \omega_+'\}$ and $\{\omega_-'', \omega_+''\}$ from the set $\Omega_1^\mp$ the lowest is the one with the smaller number from the set $\Omega_1^-$, and if these numbers are equal, the pair with the smaller number from the set $\Omega_1^+$. The order introduced in this way on the set $\Omega_2$ is linear and therefore the numbering $j_2$ is uniquely extended to entire the set $\Omega_2$.

Let us now denote by $\rho_1$ the order-preserving isomorphic mapping of the set $\widehat{\Omega_1^0 \cup \Omega_1^-}$ onto the set $\Omega_1^0 \cup \Omega_1^-$ and, using the Kronecker symbol $\delta$, we define the function $\mathcal{G}_1$ on the product $\Omega_1 \times \Omega_2$ and the linear operator $G_1: R^{\Omega_1} \to R^{\Omega_2}$ by the relations

$$\mathcal{G}_1(\omega_2, \omega_1) \stackrel{\text{def}}{=} \begin{cases} \delta(\omega_1, \rho_1(\omega_2)), & \omega_2 \in \widehat{\Omega_1^0 \cup \Omega_1^-} \\ g^1(\omega_1^+)\delta(\omega_1, \omega_1^-) - g^1(\omega_1^-)\delta(\omega_1, \omega_1^+), & \omega_2 = \{\omega_1^-, \omega_1^+\} \in \Omega_1^\mp \end{cases} \quad (2)$$

$$[G_1 h](\omega_2) \stackrel{\text{def}}{=} \sum_{\omega_1} \mathcal{G}_1(\omega_2, \omega_1) h(\omega_1) \quad h \in R^{\Omega_1} \quad (3)$$

Let us define the function $g^2$ on the set $\Omega_2$, setting

$$g^2(\omega_2) \stackrel{\text{def}}{=} [G_1 f^2](\omega_2) = \begin{cases} f^2(\rho_1(\omega_2)), & \omega_2 \in \widehat{\Omega_1^0 \cup \Omega_1^-} \\ g^1(\omega_1^+)f^2(\omega_1^-) - g^1(\omega_1^-)f^2(\omega_1^+), & \omega_2 = \{\omega_1^-, \omega_1^+\}, \end{cases} \quad (4)$$

Using the function $g^2$ instead of $g^1$ and replacing index 1 with index 2, we can define sets

$\Omega_2^0$, $\Omega_2^-$, $\Omega_2^+$, $\Omega_2^\mp$, $\widehat{\Omega_2^0 \cup \Omega_2^-}$, $\Omega_3 \stackrel{\text{def}}{=} \widehat{\Omega_2^0 \cup \Omega_2^-} \cup \Omega_2^\mp$, as well as the linear order $j_3$ on the set $\Omega_3$, the order-preserving isomorphism $\rho_2: \widehat{\Omega_2^0 \cup \Omega_2^-} \to \Omega_2^0 \cup \Omega_2^-$ and the function $\mathcal{G}_2$ on the product $\Omega_2 \times \Omega_3$:

$$\mathcal{G}_2(\omega_3, \omega_2) \stackrel{\text{def}}{=} \begin{cases} \delta(\omega_2, \rho_2(\omega_3)), & \omega_3 \in \widehat{\Omega_2^0 \cup \Omega_2^-} \\ g^2(\omega_2^+)\delta(\omega_2, \omega_2^-) - g^2(\omega_2^-)\delta(\omega_2, \omega_2^+), & \omega_3 = \{\omega_2^-, \omega_2^+\} \end{cases} \quad (5)$$

and then define the linear operator $G_2: R^{\Omega_2} \to R^{\Omega_3}$, $[G_2 h](\omega_3) \stackrel{\text{def}}{=} \sum_{\omega_2} \mathcal{G}_2(\omega_3, \omega_2) h(\omega_2)$ and the function $g^3 \stackrel{\text{def}}{=} G_2 G_1 f^3$.

Acting in this way, we can construct for $l \leq m-1$ the sets $\Omega_l^0$, $\Omega_l^-$, $\Omega_l^+$, $\Omega_l^\mp = \Omega_l^- \times \Omega_l^+$, $\Omega_{l+1} \stackrel{\text{def}}{=} \widehat{\Omega_l^0 \cup \Omega_l^-} \cup \Omega_l^\mp$, numbering $j_{l+1}$ on the set $\Omega_{l+1}$, isomorphism $\rho_l$ and function $\mathcal{G}_l$ on the product $\Omega_l \times \Omega_{l+1}$:

$$\mathcal{G}_l(\omega_{l+1}, \omega_l) \stackrel{\text{def}}{=} \begin{cases} \delta(\omega_l, \rho_l(\omega_{l+1})), & \omega_{l+1} \in \widehat{\Omega_l^0 \cup \Omega_l^-} \\ g^l(\omega_l^+)\delta(\omega_l, \omega_l^-) - g^l(\omega_l^-)\delta(\omega_l, \omega_l^+), & \omega_{l+1} = \{\omega_l^-, \omega_l^+\} \end{cases} \quad (6)$$

and construct $G_l$ operators with $\mathcal{G}_l$ kernel, and also define:

operators $G(l) \stackrel{\text{def}}{=} G_l \ldots G_1$ and functions $g^l \stackrel{\text{def}}{=} G(l) f^l$ for natural numbers $l \leq m-1$.

Let us introduce the functions $e^1, e^2$ using the relations

$$e^1 \stackrel{\text{def}}{=} c(1)\delta(n_1, 1) + \ldots + c(n)\delta(n_1, n), \tag{7}$$

$$e^2 \stackrel{\text{def}}{=} \frac{1}{b_1}(a_1(1)\delta(n_1, 1) + \ldots + a_1(n)\delta(n_1, n) + \delta(n_1, n+1))$$

and define the functions $d^k$, k=1, 2 on the set $\Omega_m$, setting $d^k(\omega_m) = [G(m-1)e^k](\omega_m)$. Let us denote $G^t(m-1)$ the operator transposed to the operator $G(m-1)$, and the kernel of this operator we denote $\mathcal{G}^t(m-1|\omega_1, \omega_m)$.

**H. Uzawa's theorem [4].** Let $\Omega_m^+ \stackrel{\text{def}}{=} \{\omega_m \in \Omega_m: d^2(\omega_m) > 0\}$, $\bar{\lambda} = \max_{\omega_m \in \Omega_m^+}\{\frac{d^1(\omega_m)}{d^2(\omega_m)}\}$. Then:

1. For the function $w \in R^{\Omega_m}$, $w \geq 0$, the restriction of the function $G^t(m-1)w$ to the set $\{1, \ldots, n\}$ is a solution to the linear programming problem under consideration if and only if $\sum_{\omega_m \in \Omega_m^+} w(\omega_m) d^2(\omega_m) = 1$ and $\sum_{\omega_m \in \Omega_m^+} w(\omega_m) d^2(\omega_m) = 1$ for $\omega_m \notin \Omega_m^+$.

2. For each value of the parameter $\omega_m$ from the set $n\Omega_m$, there is one and only one vertex of the polyhedron $V_m$, defined by the constraints of the problem under consideration, and the Cartesian coordinates of this vertex are specified by a vector with components $x(\omega_m|j)$:

$$x(\omega_m|j) \stackrel{\text{def}}{=} 1/d^2(\omega_m)\, G^t(m-1|\omega_m, j), j = 1, \ldots, n. \tag{8}$$

From this theorem it follows that to obtain all solutions to the linear programming problem under consideration, it is sufficient to construct an algorithm for calculating the function $d^1(\omega_m)/d^2(\omega_m)$ and construct functions of the argument $\omega_1$ obtained by substituting into the function $G^t(m-1|\omega_m, \omega_1)$ values of the variable $\omega_m = \bar{\omega}_m$ at which the ratio $d^1(\bar{\omega}_m)/d^2(\bar{\omega}_m)$ is maximum. The set of all solutions to the problem under consideration is the convex hull of the set of functions mentioned above.

## Construction of solutions.

**Theorem.** There are relations for any natural $l$, $m - 1 > l > 1$

$$\mathcal{G}(l|\omega_{l+1}, \omega_1) = \mathcal{G}(l-1|\rho_l(\omega_{l+1}), \omega_1), \text{ if } \omega_{l+1} \in \widehat{\Omega_l^0 \cup \Omega_l^-} \tag{9}$$

$$\mathcal{G}(l|\omega_{l+1}, \omega_1) = g^l(\omega_l^+)\mathcal{G}(l-1|\omega_l^-, \omega_1) - g^l(\omega_l^-)\mathcal{G}(l-1|\omega_l^+, \omega_1) \text{ if } \omega_l=\{\omega_l^-, \omega_l^+\} \tag{10}$$

and for an arbitrary function $h$ on the set $\Omega_1$ and the equalities hold

$$[G(l)h](\omega_{l+1}) = [G(l-1)h](\rho_l(\omega_{l+1}), \text{ if } \omega_{l+1} \in \widehat{\Omega_l^0 \cup \Omega_l^-} \tag{11}$$

$$[G(l)h](\{\omega_l^-, \omega_l^+\}) = g^l(\omega_l^+)[G(l-1)h](\omega_l^-) - g^l(\omega_l^-)[G(l-1)h](\omega_l^+) \tag{12}$$

**Proof.** Relations (9)-(12) follow from relation (6) and the equality $G(l) = G_l G(l-1)$.

We get the following algorithms using H. Uzawa's theorem.

**Algorithm for calculating the functions $g^l$, $d^1(\omega_m)$, $d^2(\omega_m)$.**

**Calculation of functions $g^l$, $l = 1, \ldots, m$ for $m > 1$.**

Let us set $f_{(s)}^l = G(s-1)f^l$ for $m - 1 \geq l \geq 2$ и $2 \leq s \leq l$. Then $g^l = f_{(l-1)}^l$.

From relations (11), (12) for the function $\mathcal{G}(l|\omega_{l+1}, \omega_1)$ recurrence relations follow

$$f^l_{(s)}(\omega_s) = f^{l-1}_{(s)}(\rho_s(\omega_s)), \text{ if } \omega_s \in \widehat{\Omega^0_s \cup \Omega^-_s}$$

$$f^l_{(s)}(\omega_s) = g^{s-1}(\omega^+_{s-1})f^l_{(s-1)}(\omega^-_{s-1}) - g^s(\omega^-_{s-1})f^l_{(s-1)}(\omega^+_{s-1}) \text{ if } \omega_s = \{\omega^-_{s-1}, \omega^+_{s-1}\}$$

giving an algorithm for calculating functions $g^l$.

For a function $h$ on the set $\Omega_1$ we denote $h_{(l+1)} \stackrel{\text{def}}{=} G(l)h$, $l = 1, \ldots, m-1$ and $h_{(1)} \stackrel{\text{def}}{=} h$.

$$h_{(l+1)}(\omega_{l+1}) = h_{(l)}(\rho_l(\omega_{l+1})) \text{ if } \omega_{l+1} \in \widehat{\Omega^0_l \cup \Omega^-_l} \qquad (13)$$

$$h_{(l+1)}(\{\omega^-_l, \omega^+_l\}) = g^l(\omega^+_l)h_{(l)}(\omega^-_l) - g^l(\omega^-_l)h_{(l)}(\omega^+_l), \text{ if } \{\omega^-_l, \omega^+_l\} \in \Omega^{\mp}_l \qquad (14)$$

Assuming $h = e^1$, $h = e^2$, we obtain an algorithm for calculating the functions $d^1(\omega_m)$, $d^2(\omega_m)$ respectively.

**Algorithm for calculating the coordinates of the vertices of the polyhedron $V_m$.**

The coordinates $x_j(\omega_m), j = 1, \ldots, n$ of the vertices of the polyhedron $V_m$ generated by the constraints can be calculated using the formula

$$x_j(\omega_l) \stackrel{\text{def}}{=} 1/d^2(\omega_l) \, G^t(l-1|\omega_l, j), j = 1, \ldots, n+1; l = 2, \ldots, m$$

We have following two relations from the equality (2):

$$x_j(k) = 1/d^2(k) \, \mathcal{G}(1|\rho_1(k), j) = x_k(\rho_1(k))/d^2(k) \text{ if } k \in \widehat{\Omega^0_1 \cup \Omega^-_1}$$

$$x_j(\omega_2) = 1/d^2(\omega_2)[f^1(\omega^+_1)x_j(\omega^-_1) - f^1(\omega^-_1)x_j(\omega^+_1)], \text{ if } \omega_2 = \{\omega^-_1, \omega^+_1\},$$

So we get the opportunity to calculate the coordinates $x_j(\omega_{l+1})$ using the relations (11), (12) in $m-1$ steps:

$$x_j(\omega_{l+1}) = 1/d^2(\omega_{l+1})x_j(\rho_l(\omega_{l+1})), \text{ if } \omega_{l+1} \in \widehat{\Omega^0_l \cup \Omega^-_l} \qquad (15)$$

$$x_j(\omega_{l+1}) = 1/d^2(\omega_{l+1})[g^l(\omega^+_l)x_j(\omega^-_l) - g^l(\omega^-_l)x_j(\omega^+_l)], \text{ if } \omega_{l+1} = \{\omega^-_l, \omega^+_l\} \qquad (16)$$

As a result, with $l = m-1$ we obtain the mentioned coordinates.

**Notation.** Let us now introduce indexing of elements from the sets $\Omega_l$:

for natural numbers $l, n_1$, we denote $[n_1]_{l+1}$ the element $\omega_{l+1} = n_1$ from the set $\widehat{\Omega^0_l \cup \Omega^-_l}$ and we denote $[n^-, n^+]_{l+1}$ element $\{\omega^-_l, \omega^+_l\}$ from the set $\Omega^{\mp}_l$ if $j_l(\omega^-_l) = n^-, j_l(\omega^+_l) = n^+$.

### 3. J.D. Wiliams' example [5].

Consider a zero-sum game defined by a payoff table:

| 4 | 3 | 3 | 2 | 2 | 6 |
|---|---|---|---|---|---|
| 0 | 7 | 3 | 6 | 2 | 2 |
| 6 | 0 | 4 | 2 | 6 | 2 |

The initial data of the task has the form:

$e^1 = \delta(\omega_1, 1) + \delta(\omega_1, 2) + \delta(\omega_1, 3) + \delta(\omega_1, 4) + \delta(\omega_1, 5) + \delta(\omega_1, 6)$

$e^2 = 4\delta(\omega_1, 1) + 3\delta(\omega_1, 2) + 3\delta(\omega_1, 3) + 2\delta(\omega_1, 4) + 2\delta(\omega_1, 5) + 6\delta(\omega_1, 6) + \delta(\omega_1, 7)$

$f^1 = -4\delta(\omega_1, 1) + 4\delta(\omega_1, 2) + 4\delta(\omega_1, 4) - 4\delta(\omega_1, 6) - \delta(\omega_1, 7)$

$f^2 = =2\delta(\omega_1, 1) - 3\delta(\omega_1, 2) + \delta(\omega_1, 3) + 4\delta(\omega_1, 5) - 4\delta(\omega_1, 6) - \delta(\omega_1, 7)$

and the index sets have the form:

$\Omega_1^0 \cup \Omega_1^- = \{1,3,5,6,7\}, \Omega_1^- = \{1,6,7\}, \Omega_1^+ = \{2,4\}, \widehat{\Omega_1^0 \cup \Omega_1^-} = \{[1]_2, [2]_2, [3]_2, [4]_2, [5]_2\}$

$\Omega_1^- \times \Omega_1^+ = \{[1,2]_2, [1,4]_2, [6,2]_2, [6,4]_2, [7,2]_2, [7,4]_2\}$

$\Omega_2 = \{[1]_2, [2]_2, [3]_2, [4]_2, [5]_2, [1,2]_2, [1,4]_2, [6,2]_2, [6,4]_2, [7,2]_2, [7,4]_2\}$

and function $\rho_1$: given by the table

| $\omega_2$ | $[1]_2$ | $[2]_2$ | $[3]_2$ | $[4]_2$ | $[5]_2$ |
|---|---|---|---|---|---|
| $\rho_1(\omega_2)$ | 1 | 3 | 5 | 6 | 7 |

Applying relations (9), (10) we will ultimately obtain a table of the values of the variable $\omega_2$ and the values of the functions $g^2, [G(1)e^1](\omega_2), [G(1)e^2](\omega_2)$

**Table 1.**

| $\omega_2$ | $[1]_2$ | $[2]_2$ | $[3]_2$ | $[4]_2$ | $[5]_2$ | $[1,2]_2$ | $[1,4]_2$ | $[6,2]_2$ | $[6,4]_2$ | $[7,2]_2$ | $[7,4]_2$ |
|---|---|---|---|---|---|---|---|---|---|---|---|
| $g^2$ | 2 | 1 | 4 | -4 | -1 | -4 | 8 | -28 | -16 | -7 | -4 |
| $G(1)e^1$ | 1 | 1 | 1 | 1 | 0 | 8 | 8 | 8 | 8 | 1 | 1 |
| $G(1)e^2$ | 4 | 3 | 2 | 6 | 1 | 28 | 24 | 36 | 32 | 7 | 6 |

Then using the function $g^2(\omega_2)$, we get a table to display $\rho_2$

| $\omega_3$ | $[1]_3$ | $[2]_3$ | $[3]_3$ | $[4]_3$ | $[5]_3$ | $[6]_3$ | $[7]_3$ |
|---|---|---|---|---|---|---|---|
| $\rho_2(\omega_3)$ | $[4]_2$ | $[5]_2$ | $[1,2]_2$ | $[6,2]_2$ | $[6,4]_2$ | $[7,2]_2$ | $[7,4]_2$ |

Using the table for restrictions of functions $d^k$ on the set $\widehat{\Omega_2^0 \cup \Omega_2^-}$

| | $[1]_3$ | $[2]_3$ | $[3]_3$ | $[4]_3$ | $[5]_3$ | $[6]_3$ | $[7]_3$ |
|---|---|---|---|---|---|---|---|
| $d^1 = G(2)e^1$ | 1 | 0 | 8 | 8 | 8 | 1 | 1 |
| $d^2 = G(2)e^2$ | 6 | 1 | 28 | 36 | 32 | 7 | 6 |
| $d^1/d^2$ | 1/6 | 0 | 2/7 | 2/9 | 1/4 | 1/7 | 1/6 |

we obtain the index values for 7 more vertices of the polyhedron $V_3$ generated by the restrictions

| $[1]_3$ | $[2]_3$ | $[3]_3$ | $[4]_3$ | $[5]_3$ | $[6]_3$ | $[7]_3$ |
|---|---|---|---|---|---|---|

From the table for the function $g^2$ we obtain the equalities

$$\Omega_2^- = \{[4]_2, [5]_2, [\,1,2]_2, [\,6,2]_2, [\,6,4]_2, [\,7,2]_2, [\,7,4]_2\}, \quad \Omega_2^+ = \{[1]_2, [2]_2, [3]_2, [\,1,4]_2\}$$

and the table of indices for the remaining 28 vertices of the polyhedron $V_3$

|  | $[4]_2$ | $[5]_2$ | $[\,1,2]_2$ | $[\,6,2]_2$ | $[\,6,4]_2$ | $[\,7,2]_2$ | $[\,7,4]_2$ |
|---|---|---|---|---|---|---|---|
| $[1]_2$ | $[[4]_2,[1]_2]_3$ | $[[5]_2,[1]_2]_3$ | $[[1,2]_2,[1]_2]_3$ | $[[6,2]_2,[1]_2]_3$ | $[[6,4]_2,[1]_2]_3$ | $[[7,2]_2,[1]_2]_3$ | $[[7,4]_2,[1]_2]_3$ |
| $[2]_2$ | $[[4]_2,[2]_2]_3$ | $[[5]_2,[2]_2]_3$ | $[[1,2]_2,[2]_2]_3$ | $[[6,2]_2,[2]_2]_3$ | $[[6,4]_2,[2]_2]_3$ | $[[7,2]_2,[2]_2]_3$ | $[[7,4]_2,[2]_2]_3$ |
| $[3]_2$ | $[[4]_2,[3]_2]_3$ | $[[5]_2,[3]_2]_3$ | $[[1,2]_2,[3]_2]_3$ | $[[6,2]_2,[3]_2]_3$ | $[[6,4]_2,[3]_2]_3$ | $[[7,2]_2,[3]_2]_3$ | $[[7,4]_2,[3]_2]_3$ |
| $[\,1,4]_2$ | $[[4]_2,[1,4]_2]_3$ | $[[5]_2,[1,4]_2]_3$ | $[[1,2]_2,[1,4]_2]_3$ | $[[6,2]_2,[1,4]_2]_3$ | $[[6,4]_2,[1,4]_2]_3$ | $[[7,2]_2,[1,4]_2]_3$ | $[[7,4,[1,4]_2]_3$ |

Thus we get that the polyhedron $V_3$ has 35 vertices. To calculate the values of the functions $d^1, d^2$ for the argument $\omega_3$ from the set $\Omega_3^\mp$, we use relation (12) as well as **Table 1.**

| $\omega_2$ | $[1]_2$ | $[2]_2$ | $[3]_2$ | $[4]_2$ | $[5]_2$ | $[\,1,2]_2$ | $[\,1,4]_2$ | $[6,2]_2$ | $[6,4]_2$ | $[7,2]_2$ | $[7,4]_2$ |
|---|---|---|---|---|---|---|---|---|---|---|---|
| $g^2$ | 2 | 1 | 4 | -4 | -1 | -4 | 8 | -28 | -16 | -7 | -4 |
| $h^2 = G(1)e^2$ | 4 | 3 | 2 | 6 | 1 | 28 | 24 | 36 | 32 | 7 | 6 |

As a result, we get the following table

|  | $\omega_3 = \{\omega_2^- \times \omega_2^+\}$ | $d^1$ | $d^2$ | $d^1/d^2$ |  | $\omega_3 = \{\omega_2^- \times \omega_2^+\}$ | $d^1$ | $d^2$ | $d^1/d^2$ |
|---|---|---|---|---|---|---|---|---|---|
| 1 | $[[4]_2,[1]_2]_3$ | 6 | 28 | 3/14 | 15 | $[[6,2]_2,[3]_2]_3$ | 60 | 200 | 3/10 |
| 2 | $[[4]_2,[2]_2]_3$ | 5 | 18 | 5/18 | 16 | $[[6,2]_2,[1,4]_2]_3$ | 288 | 960 | 3/10 |
| 3 | $[[4]_2,[3]_2]_3$ | 8 | 32 | 1/4 | 17 | $[[6,4]_2,[1]_2]_3$ | 32 | 128 | 1/4 |
| 4 | $[[4]_2,[1,4]_2]_3$ | 40 | 144 | 5/18 | 18 | $[[6,4]_2,[2]_2]_3$ | 24 | 80 | 3/10 |
| 5 | $[[5]_2,[1]_2]_3$ | 1 | 6 | 1/6 | 19 | $[[6,4]_2,[3]_2]_3$ | 48 | 160 | 3/10 |
| 6 | $[[5]_2,[2]_2]_3$ | 1 | 4 | 1/4 | 20 | $[[6,4]_2,[1,4]_2]_3$ | 192 | 640 | 3/10 |
| 7 | $[[5]_2,[3]_2]_3$ | 1 | 6 | 1/6 | 21 | $[[7,2]_2,[1]_2]_3$ | 9 | 42 | 3/14 |
| 8 | $[[5]_2,[1,4]_2]_3$ | 8 | 32 | 1/4 | 22 | $[[7,2]_2,[2]_2]_3$ | 8 | 28 | 2/7 |
| 9 | $[[1,2]_2,[1]_2]_3$ | 8 | 72 | 1/9 | 23 | $[[7,2]_2,[3]_2]_3$ | 11 | 42 | 11/42 |
| 10 | $[[1,2]_2,[2]_2]_3$ | 12 | 40 | 3/10 | 24 | $[[7,2]_2,[1,4]_2]_3$ | 64 | 224 | 2/7 |
| 11 | $[[1,2]_2,[3]_2]_3$ | 36 | 120 | 3/10 | 25 | $[[7,4]_2,[1]_2]_3$ | 6 | 28 | 3/14 |
| 12 | $[[1,2]_2,[1,4]_2]_3$ | 96 | 320 | 3/10 | 26 | $[[7,4]_2,[2]_2]_3$ | 5 | 18 | 5/18 |
| 13 | $[[6,2]_2,[1]_2]_3$ | 44 | 184 | 11/46 | 27 | $[[7,4]_2,[3]_2]_3$ | 8 | 32 | 1/4 |
| 14 | $[[6,2]_2,[2]_2]_3$ | 36 | 120 | 3/10 | 28 | $[[7,4,[1,4]_2]_3$ | 40 | 144 | 5/18 |

From this table it can be seen that:

1. the price of the game is 3/10,

2. optimal functions correspond to the following vertices:

$[[\,1,2]_2,[2]_2]_3, [[\,1,2]_2,[3]_2]_3, [[\,6,2]_2,[1]_2]_3, [[\,6,2]_2,[2]_2]_3, [[\,6,2]_2,[3]_2]_3, [[\,6,4]_2,[1]_2]_3,$
$[[\,6,4]_2,[2]_2]_3, [[\,6,4]_2,[3]_2]_3, [[\,6,4]_2,[1.4]_2]_3.$

As a result, using algorithm (12) (13) we obtain a table for nine extremal functions.

|   |   |   | $d^1$ | $d^2$ | $d^1/d^2$ |
|---|---|---|---|---|---|
| 1 | $[[1,2]_2, [2]_2]_3$ | $4\delta(n_1, 1) + 4\delta(n_1, 2) + 4\delta(n_1, 3)$ | 12 | 40 | 3/10 |
| 2 | $[[1,2]_2, [3]_2]_3$ | $16\delta(n_1, 1) + 4\delta(n_1, 2) + 4\delta(n_1, 5)$ | 36 | 120 | 3/10 |
| 3 | $[[1,2]_2, [1,4]_2]_3$ | $48\delta(\omega_1, 1) + 32\delta(\omega_1, 2) + 16\delta(\omega_1, 4)$ | 96 | 320 | 3/10 |
| 4 | $[[6,2]_2, [2]_2]_3$ | $4\delta(\omega_1, 2) + 28\delta(\omega_1, 3) + 4\delta(\omega_1, 6)$ | 36 | 120 | 3/10 |
| 5 | $[[6,2]_2, [3]_2]_3$ | $16\delta(\omega_1, 2) + 28\delta(\omega_1, 5) + 16\delta(\omega_1, 6)$ | 60 | 200 | 3/10 |
| 6 | $[[6,2]_2, [1,4]_2]_3$ | $112\delta(\omega_1, 1) + 32\delta(\omega_1, 2) + 112\delta(\omega_1, 4) + 32\delta(\omega_1, 6)$ | 288 | 960 | 3/10 |
| 7 | $[[6,4]_2, [2]_2]_3$ | $16\delta(\omega_1, 3) + 4\delta(\omega_1, 4) + 4\delta(\omega_1, 6)$ | 24 | 80 | 3/10 |
| 8 | $[[6,4]_2, [3]_2]_3$ | $16\delta(n_1, 4) + 16\delta(n_1, 5) + 16\delta(n_1, 6)$ | 48 | 160 | 3/10 |
| 9 | $[[6,4]_2, [1,4]_2]_3$ | $64\delta(\omega_1, 1) + 48\delta(\omega_1, 4) + 16\delta(\omega_1, 6))$ | 12 | 40 | 3/10 |

The table of nine optimal strategies of the second player in vector form looks like

| 1 | 2 | 3 | 4 | 5 | 6 | 7 | 8 | 9 |
|---|---|---|---|---|---|---|---|---|
| 1/3 | 4/9 | 1/2 | 0 | 0 | 7/18 | 0 | 0 | 1/2 |
| 1/3 | 4/9 | 1/3 | 1/9 | 4/15 | 1/9 | 0 | 0 | 0 |
| 1/3 | 0 | 0 | 7/9 | 0 | 0 | 2/3 | 0 | 0 |
| 0 | 0 | 1/6 | 0 | 0 | 7/18 | 1/6 | 1/3 | 3/8 |
| 0 | 1/9 | 0 | 0 | 7/15 | 0 | 0 | 1/3 | 0 |
| 0 | 0 | 0 | 1/9 | 4/15 | 1/9 | 1/6 | 1/3 | 1/8 |

Using algorithm (15), (16) one can also calculate the coordinates of all other vertices of the polyhedron $V_3$.

**References.**